\documentclass{article}
\usepackage{amssymb, latexsym, amsmath, eucal, graphicx, epsfig, amsthm,bm, fancyhdr}
\usepackage[top=1in, bottom=1in, left=.75in, right=.75in]{geometry}
\newtheorem{theorem}{Theorem}
\newtheorem{lemma}{Lemma}

\newtheorem{corollary}{Corollary}

\everymath{\displaystyle}
\title{Updating the Inverse of a Matrix When Removing the $i$th Row and Column with an Application to Disease Modeling}
\author{Cody Palmer}

\begin{document}

\maketitle

\begin{abstract}
	The Sherman-Woodbury-Morrison (SWM) formula gives an explicit formula for the inverse perturbation of a matrix
	in terms of the inverse of the original matrix and the perturbation.  This formula is useful for numerical
	applications.  We have produced similar results, giving an expression for the inverse of a matrix when the $i$th
	row and column are removed.  However, our expression involves taking a limit, which inhibits use in similar
	applications as the SWM formula.  However, using our expression to find an analytical result on the spectral
	radius of a special product of two matrices leads to an application. In particular, we find a way to compute
	the fundamental reproductive ratio of a relapsing disease being spread by a vector among two species of host that
	undergo a different number of relapses.
\end{abstract}

\section{Introduction}
The process of
computing the inverse of matrix after altering it is known as updating the inverse of a matrix \cite{hager1989updating}.
The most famous example of this process are the Sherman-Morrison-Woodbury formula which gives a closed form expression for
the inverse of a perturbation of a matrix in terms of its original inverse.  In this paper we will be representing the
inverse of a matrix when the $i$th row and column are removed as a limit involving the original inverse.
While such methods have numerical applications \cite{hager1989updating},
our method, since it contains a limit, is going to have a more analytical usage.  In particular, after some
introduction, we will show how this
result can be used to compute the fundamental reproductive ratio for a model involving a relapsing disease being spread 
among two host species by a vector.

\section{Main Results}
We begin by investigating the determinant of a square matrix $A$ as a diagonal element tends to 
$\infty$. Note that $A_{[i,j]}$ represents the matrix formed by removing the $i$th row and $j$th column, and
we will denote a particular element of a matrix with a lower case letter corresponding to the matrix e.g. $a_{mn}$ is the
$m,n$th element of $A$.  We will also occasionally use $(B)_{ij}$ to represent the $i,j$th element of $B$.
\begin{lemma}\label{detlemma}
	Let $A$ be an $n \times n$ matrix and suppose that $A_{[i,i]}$ is nonsingular.  Then $\lim_{a_{ii} \to \infty} \det{A} = \pm \infty$.
\end{lemma}
\begin{proof} By Proposition 2.7.5 of \cite{matmat}:
\[
	\det{A} = \sum_{k=1}^n (-1)^{i+k}a_{ik}\det(A_{[i,k]}) = a_{ii}\det{A_{[i,i]}} +
 \sum_{k \neq i} (-1)^{i+k}a_{ik}\det(A_{[i,k]})
\]
The last sum does not involve $a_{ii}$ and thus has a fixed value as $a_{ii} \to \infty$.  Since $A_{[i,i]}$ is nonsingular
it has a nonzero determinant, and thus the leading term of the previous sum goes to $\pm\infty$ depending on the sign of
$\det{A_{[i,i]}}$.
\end{proof}

As a result of this lemma we can see that there is a sufficiently large value of $a_{ii}$ that makes $A$ invertible, and the
matrix remains invertible for all further values.  Hence, the hypothesis of $A$ being invertible is not needed in the 
next result, which tells us how to construct the inverse of $A_{[i,i]}$ from $A^{-1}$:
\begin{theorem} \label{limlemma} If $A_{[i,i]}$ is nonsingular then
\[
	(A_{[i,i]})^{-1} = \lim_{a_{ii} \to \infty} (A^{-1})_{[i,i]}
\]
Furthermore,
\[
	\lim_{a_{ii} \to \infty} (A^{-1})_{ik} = \lim_{a_{ii} \to \infty} (A^{-1})_{ki} = 0
\]
\end{theorem}
\begin{proof}
	We will need to consider this proof in four cases.   The proof technique in each case is the same, though the indexing in each is different.
	Throughout let $B^{jk} = (b_{pq}) = A_{[j,k]}$ for 
	$1 \leq p,q \leq n-1$.  We will repeatedly use Corollary 2.7.6 of \cite{matmat} which is a formula for the $ij$ element of the inverse of a matrix.
	Also, the use of a ``$\ldots$" denotes terms of a sum that do not involve $a_{ii}$.

	\textbf{Case 1.}  Assume that $1 \leq j,k < i$.
	Then, on the one hand
	\[
		\lim_{a_{ii} \to \infty} (A^{-1})_{kj} = \lim_{a_{ii} \to \infty} \frac{(-1)^{k+j} \det B^{jk}}{\det A}
		\]
		\[
			=\lim_{a_{ii} \to \infty} \frac{(-1)^{k+j} \sum_{l=1}^{n-1} (-1)^{i-1+l} b_{i-1,l} 
			\det B^{jk}_{[i-1,l]}}{a_{ii}\det A_{[i,i]} + \ldots}
	\]
	after expanding $\det B^{jk}$ along its $i-1$ row.  Now we want to identify the term that has $a_{ii}$ in it.  Note that
	since $j,k <i$ we have that $a_{ii} = b_{i-1,i-1}$.  So then we let $l=i-1$ and we have
	\[
		\lim_{a_{ii} \to \infty} (A^{-1})_{kj} = \lim_{a_{ii} \to \infty} 
		\frac{(-1)^{k+j} (-1)^{2i-2} a_{ii} \det B^{jk}_{[i-1,i-1]} + \ldots}
		{ a_{ii} \det A_{[i,i]} + \ldots}	
	\]
	\[
		= \frac{(-1)^{k+j}   \det B^{jk}_{[i-1,i-1]} }
		{\det A_{[i,i]}}
	\]
	On the other hand
	\[
		((A_{[i,i]})^{-1})_{kj} = \frac{(-1)^{k+j} \det (A_{[i,i]})_{[j,k]}}{\det A_{[i,i]}}
	\]
	Since $j,k < i$ we have that $(A_{[i,i]})_{[j,k]} = (A_{[j,k]})_{[i-1,i-1]} = B^{jk}_{[i-1,i-1]}$ so that
	\[
		((A_{[i,i]})^{-1})_{kj} = \frac{(-1)^{k+j} \det B^{jk}_{[i-1,i-1]}}{\det A_{[i,i]}}
	\]
	Thus when $j,k<i$ we have
	\[
		((A_{[i,i]})^{-1})_{kj} = \lim_{a_{ii} \to \infty} (A^{-1})_{kj}
	\]

	\textbf{Case 2.} $n-1 \geq k,j \geq i$.  On the one hand
	\[
		\lim_{a_{ii} \to \infty} (A^{-1})_{k+1,j+1} = \lim_{a_{ii} \to \infty} \frac{(-1)^{k+j+2} \det B^{j+1,k+1}}
		{\det A}
	\]
	\[
		= \lim_{a_{ii} \to \infty} \frac{(-1)^{k+j} \sum_{l=1}^{n-1} (-1)^{i+l} b_{i,l} 
			\det B^{j+1,k+1}_{[i,l]}}{a_{ii}\det A_{[i,i]} + \ldots}
	\]
	after expanding $\det B^{j+1,k+1}$ along its $i$th row.  We have that $a_{ii} = b_{ii}$.  
	So then we let $l=i$ and we have
	\[
		\lim_{a_{ii} \to \infty} (A^{-1})_{kj} = \lim_{a_{ii} \to \infty} 
		\frac{(-1)^{k+j} (-1)^{2i} a_{ii} \det B^{j+1,k+1}_{[i,i]} + \ldots}
		{ a_{ii} \det A_{[i,i]} + \ldots}	
	\]
	\[
		= \frac{(-1)^{k+j}   \det B^{j+1,k+1}_{[i,i]} }
		{\det A_{[i,i]}}
	\]
	On the other hand
	\[
		((A_{[i,i]})^{-1})_{kj} = \frac{(-1)^{k+j} \det (A_{[i,i]})_{[j,k]}}{\det A_{[i,i]}}
	\]
	Since $j,k \geq i$ we have that $(A_{[i,i]})_{[j,k]} = (A_{[j+1,k+1]})_{[i,i]} = B^{j+1,k+1}_{[i,i]}$ so that
	\[
		((A_{[i,i]})^{-1})_{kj} = \frac{(-1)^{k+j} \det B^{j+1,k+1}_{[i,i]}}{\det A_{[i,i]}}
	\]
	Thus when $j,k\geq i$ we have
	\[
		((A_{[i,i]})^{-1})_{kj} = \lim_{a_{ii} \to \infty} (A^{-1})_{k+1,j+1}
	\]

	\textbf{Case 3.} $k=i, j <i$. On the one hand
	\[
		\lim_{a_{ii} \to \infty} (A^{-1})_{i+1,j} = \lim_{a_{ii} \to \infty} \frac{(-1)^{i+j+1} \det B^{j,i+1}}
		{\det A}
	\]
	\[
		= \lim_{a_{ii} \to \infty} \frac{(-1)^{i+j+1} \sum_{l=1}^{n-1} (-1)^{i-1+l} b_{i-1,l} 
			\det B^{j,i+1}_{[i-1,l]}}{a_{ii}\det A_{[i,i]} + \ldots}
	\]
	after expanding $\det B^{j,i+1}$ along its $i-1$th row. Since $j <i$  we have that $a_{ii} = b_{i-1,i}$.  
	So then we let $l=i$ and we have
	\[
		\lim_{a_{ii} \to \infty} (A^{-1})_{i+1,j} = \lim_{a_{ii} \to \infty} 
		\frac{(-1)^{i+j+1} (-1)^{2i-1} a_{ii} \det B^{j,i+1}_{[i-1,i]} + \ldots}
		{ a_{ii} \det A_{[i,i]} + \ldots}	
	\]
	\[
		= \frac{(-1)^{i+j}   \det B^{j,i+1}_{[i-1,i]} }
		{\det A_{[i,i]}}
	\]
	On the other hand
	\[
		((A_{[i,i]})^{-1})_{ij} = \frac{(-1)^{i+j} \det (A_{[i,i]})_{[j,i]}}{\det A_{[i,i]}}
	\]
	Since $j<i$ we have that $(A_{[i,i]})_{[j,i]} = (A_{[j,i+1]})_{[i-1,i]} = B^{j,i+1}_{[i-1,i]}$ so that
	\[
		((A_{[i,i]})^{-1})_{ij} = \frac{(-1)^{i+j} \det B^{j,i+1}_{[i-1,i]}}{\det A_{[i,i]}}
	\]
	Thus when $j<i$ we have
	\[
		((A_{[i,i]})^{-1})_{ij} = \lim_{a_{ii} \to \infty} (A^{-1})_{i+1,j}
	\]

	\textbf{Case 4.} $k<i, j=i$. On the one hand
	\[
		\lim_{a_{ii} \to \infty} (A^{-1})_{k,i+1} = \lim_{a_{ii} \to \infty} \frac{(-1)^{i+k+1} \det B^{i+1,k}}
		{\det A}
	\]
	\[
		= \lim_{a_{ii} \to \infty} \frac{(-1)^{i+j+1} \sum_{l=1}^{n-1} (-1)^{i+l} b_{i,l} 
			\det B^{i+1,k}_{[i,l]}}{a_{ii}\det A_{[i,i]} + \ldots}
	\]
	after expanding $\det B^{i+1,k}$ along its $i$th row. Since $k <i$  we have that $a_{ii} = b_{i,i-1}$.  
	So then we let $l=i-1$ and we have
	\[
		\lim_{a_{ii} \to \infty} (A^{-1})_{k,i+1} = \lim_{a_{ii} \to \infty} 
		\frac{(-1)^{i+k+1} (-1)^{2i-1} a_{ii} \det B^{i+1,k}_{[i,i-1]} + \ldots}
		{ a_{ii} \det A_{[i,i]} + \ldots}	
	\]
	\[
		= \frac{(-1)^{i+k}   \det B^{i+1,k}_{[i,i-1]} }
		{\det A_{[i,i]}}
	\]
	On the other hand
	\[
		((A_{[i,i]})^{-1})_{ki} = \frac{(-1)^{i+k} \det (A_{[i,i]})_{[i,k]}}{\det A_{[i,i]}}
	\]
	Since $k<i$ we have that $(A_{[i,i]})_{[i,k]} = (A_{[i+1,k]})_{[i,i-1]} = B^{i+1,k}_{[i,i-1]}$ so that
	\[
		((A_{[i,i]})^{-1})_{ki} = \frac{(-1)^{i+k} \det B^{i+1,k}_{[i,i-1]}}{\det A_{[i,i]}}
	\]
	Thus when $k<i$ we have
	\[
		((A_{[i,i]})^{-1})_{ki} = \lim_{a_{ii} \to \infty} (A^{-1})_{k,i+1}
	\]
	The combination of these four cases gives the first result.

	For the second result, we again use Corollary 2.7.6 of \cite{matmat} to get that
	\[
		\lim_{a_{ii} \to \infty} (A^{-1})_{ik} = \lim_{a_{ii} \to \infty} \frac{(-1)^{k+i} \det B^{ki}}{\det A}
	\]
	$B^{ki}$ does not contain $a_{ii}$, and thus $\det B^{ik}$ remains constant for all values of $a_{ii}$, 
and by Lemma \ref{detlemma} we have that $\det A \to \pm \infty$.
As a result
\[
	\lim_{a_{ii} \to \infty} \frac{(-1)^{k+i} \det B^{ki}}{\det A} = 0
\]
The result for $(A^{-1})_{ki}$ is done in exactly the same way.
\end{proof}

Lemma \ref{detlemma} and Theorem \ref{limlemma}  will allow us to prove a result about the spectral radius of 
a special product of matrices
\begin{corollary} \label{conthe}
	Suppose that $V_{[i,i]}$ is nonsingular.  Then
	\[
		\lim_{v_{ii} \to \infty} \rho(FV^{-1}) = \rho(F_{[i,i]}(V_{[i,i]})^{-1})
	\]
\end{corollary}
\begin{proof}
	As before, because $V_{[i,i]}$ is nonsingular $V$ must be nonsigular for sufficiently large $v_{ii}$.
	Since eigenvalues are continuous with respect to the entries of a matrix, and the absolute value and
	maximum of a set of continuous functions is continuous, we have that
	\[
		\lim_{v_{ii} \to \infty} \rho(FV^{-1}) = \rho(F\lim_{v_{ii} \to \infty}V^{-1})
	\]
	Let 
	\[
	(V_{[i,i]})^{-1} = \begin{pmatrix} V_1 & V_2\\ V_3 & V_4 \end{pmatrix}
	\]
	where $V_1 \in \mathbb{R}^{(i-1) \times (i-1)}$, $V_2 \in \mathbb{R}^{(i-1) \times (n-i)}$
	$V_3 \in \mathbb{R}^{(n-i) \times (i-1)}$ and $V_4 \in \mathbb{R}^{(n-i) \times (n-i)}$.  Then Lemma
	\ref{limlemma} says that 
	\[
		\lim_{v_{ii} \to \infty}V^{-1} = \begin{pmatrix} V_1 & \bm{0}_{(i-1)\times 1} &V_2\\
		\bm{0}_{1 \times (i-1)} & 0 & \bm{0}_{1 \times (n-i)}\\
	V_3 & \bm{0}_{(n-i)\times 1} & V_4 \end{pmatrix}
	\]
	Let 
	\[
		F =\begin{pmatrix} F_1 & \bm{f}^{(1)}_{(i-1)\times 1} &F_2\\
		\bm{f}^{(2)}_{1 \times (i-1)} & f_{ii} & \bm{f}^{(3)}_{1 \times (n-i)}\\
	F_3 & \bm{f}^{(4)}_{(n-i)\times 1} & F_4 \end{pmatrix}
	\]
	where $F_1 \in \mathbb{R}^{(i-1) \times (i-1)}$, $F_2 \in \mathbb{R}^{(i-1) \times (n-i)}$
	$F_3 \in \mathbb{R}^{(n-i) \times (i-1)}$ and $F_4 \in \mathbb{R}^{(n-i) \times (n-i)}$.  This
	gives that
	\begin{multline*}
		F \lim_{v_{ii} \to \infty} V^{-1} \\=
		\begin{pmatrix} F_1V_1 + F_2V_3 & \bm{0}_{(i-1) \times 1} & F_1V_2 + F_2V_4\\
		\bm{f}^{(2)}_{1 \times (i-1)} V_1 + \bm{f}^{(3)}_{1 \times (n-1)}V_3 & 0 & 
		\bm{f}^{(2)}_{1 \times (i-1)} V_2 + \bm{f}^{(3)}_{1 \times (n-1)}V_4\\
		F_3V_1 + F_4V_3 & \bm{0}_{(n-i) \times 1} & F_3V_2+F_4V_4 \end{pmatrix}
	\end{multline*}
	We wish to compute the spectral radius of this matrix, so we set up the eigenvalue problem
	\begin{multline*}
		\det(F \lim_{v_{ii} \to \infty} V^{-1} - \lambda I_n)
		 \\= \det \begin{pmatrix}
			 F_1V_1 + F_2V_3 - \lambda I_{i-1}& \bm{0}_{(i-1) \times 1} & F_1V_2 + F_2V_4 	\\
		\bm{f}^{(2)}_{1 \times (i-1)} V_1 + \bm{f}^{(3)}_{1 \times (n-1)}V_3 & -\lambda & 
		\bm{f}^{(2)}_{1 \times (i-1)} V_2 + \bm{f}^{(3)}_{1 \times (n-1)}V_4\\
	F_3V_1 + F_4V_3 & \bm{0}_{(n-i) \times 1} & F_3V_2+F_4V_4 - \lambda I_{n-i} \end{pmatrix}
		 \\= -\lambda \det \begin{pmatrix}
		F_1V_1 + F_2V_3 - \lambda I_{i-1}& F_1V_2 + F_2V_4 	\\
		F_3V_1 + F_4V_3  & F_3V_2+F_4V_4 - \lambda I_{n-i} \end{pmatrix}
			  \\= -\lambda \det (F_{[i,i]}(V_{[i,i]})^{-1} - \lambda I_{n-1})
	\end{multline*}
	So the spectrum of $F \lim_{v_{ii} \to \infty} V^{-1}$ is 0 unioned with the spectrum of $F_{[i,i]}(V_{[i,i]})^{-1}$.
	Thus the spectral 
	radius of  $F \lim_{v_{ii} \to \infty} V^{-1}$ is the maximum of the eigenvalues of $F_{[i,i]}(V_{[i,i]})^{-1}$.
	That is 
	\[
		\rho(F \lim_{v_{ii} \to \infty} V^{-1}) = \rho(F_{[i,i]}(V_{[i,i]})^{-1})
	\]
	which gives the result.
	\end{proof}

\section{Applications}
To give an application we must first have a brief description of 
compartmental disease models, following their development in
\cite{compart}.  Suppose a population can be separated into $n$ homogeneous compartments and
the number of members in each compartment
will be represented by the vector $\bm{x} \in \mathbb{R}^n$ where the first $m$ compartments represent infected 
states while the remaining $n-m$ compartments
are uninfected states.  It is natural to insist that $\bm{x} \geq \bm{0}$ (inequality is taken componentwise) since we are 
dealing with populations.  Let $\bm{X}_s = \{\bm{x} \geq \bm{0}:
x_i = 0, i=1, \ldots m\}$ be the set of disease free states.  Let $\mathcal{F}_i(\bm{x})$ be the number of new infections 
in compartment $i$ (autonomy is assumed).  $\mathcal{V}_i^+(\bm{x})$ is
the rate of transfer of individuals into compartment $i$ and $\mathcal{V}_i^-(\bm{x})$ is the rate of transfer out of 
compartment $i$.  Assume that these functions are at least twice continuously differentiable.
The disease transmission model can be written as
\begin{equation} \label{compsys}
	\dot{x}_i = f_i(\bm{x}) = \mathcal{F}_i(\bm{x}) +\mathcal{V}_i^+(\bm{x}) - \mathcal{V}_i^-(\bm{x}) \qquad i=1,\ldots n
\end{equation}
Let $\mathcal{V}_i(\bm{x}) = \mathcal{V}_i^-(\bm{x}) -\mathcal{V}_i^+(\bm{x})$.  Suppose that $\bm{x}_0 \in \bm{X}_s$ is
also a fixed point of \eqref{compsys} then we call $\bm{x}_0$ a disease free equilibrium (DFE). Let $\mathcal{F}$ be the
vector valued function with the $\mathcal{F}_i$ as components, and $\mathcal{V}$ similarly defined.  Given
five conditions (A1-A5 of \cite{compart}) the Jacobians of $\mathcal{F}$ and $\mathcal{V}$ must take the form
\[
	D\mathcal{F}(\bm{x}_0) = \begin{pmatrix} F&0\\0&0 \end{pmatrix} \mbox{ and }
D\mathcal{V}(\bm{x}_0)= \begin{pmatrix} V &0 \\J_1&J_2 \end{pmatrix}
\]
Where $F$ and $V$ are $m \times m$.  Furthermore, under these conditions $V$ is nonsingular, which allows us 
to define the \textit{fundamental reproductive ratio}:
\[
	R_0 = \rho(FV^{-1})
\]
where $\rho$ is the spectral radius.
Informally, we can think of $R_0$ as being the average number of new infections produced by a single infected individual
\cite{reprat}.
With this interpretation in mind, it makes Theorem 2 of \cite{compart} expected: For $R_0<1$ the DFE $\bm{x_0}$ is stable
and for $R_0>1$ $\bm{x}_0$ is unstable.

Computing $R_0$ can be a difficult and tedious process, particularly when dealing with systems with large numbers of 
compartments. Recent work has computed $R_0$ for vector-borne diseases which relapse an arbitrary
number of times.  For a full description of these types of models see \cite{johnsonPLOSNTD}. 
In particular, the fundamental reproductive ratio was computed for two kinds of systems:
\begin{itemize}
	\item One host species undergoing $j-1$ relapses with one vector species spreading the disease.
	\item Two host species each undergoing $j-1$ relapses with one vector species spreading the disease.
\end{itemize}
In the first case, we will call the system uncoupled and
in the second we will call it a coupled system.  The equations describing the dynamics will not be reproduced here
but can be found in \cite{johnsonPLOSNTD}.

For notation, suppose that
the $i$th species is the only species in the system (an uncoupled system) and let 
$R_{0,i,j}$, $i=1,2$, be the reproductive ratio when the hosts undergo $j-1$ relapses and thus have $j$ infected
compartments. Let $F_{j,k}$ and $V_{j,k}$
be the Jacobians for the coupled system when the first host species undergoes $j-1$ relapses and the second undergoes $k-1$
relapses. Lastly, let $R_{0}^{j,k}$ be the reproductive ratio for the coupled system where the first species
undergoes $j-1$ relapses and the second species undergoes $k-1$ relapses.
We can write the reproductive ratio for the uncoupled systems in terms of the parameters for the model
\begin{equation} \label{uncoupledrnot}
	R_{0,i,j}= f \sqrt{\frac{c_ic_v\overline{S}_v}{\tilde{\mu}\overline{S}_i} \sum_{k=1}^{j} 
	\prod_{l=1}^k \frac{\alpha_{i,l-1}}{\alpha_{i,l} + \mu_{i,l}}}
\end{equation}
The coupled and uncoupled systems can then be related:
\begin{equation} \label{coupledrnot}
	R_0^{j,j} = \sqrt{(R_{0,1,j})^2+(R_{0,2,j})^2}
\end{equation}

A full description of the parameters is found in \cite{johnsonPLOSNTD}, but the relevant portion for our work here
will be to note that the first host species leaves the $j$th infected compartment at a rate $\alpha_{1,j}$ and that
$\alpha_{1,j}\to \infty$ implies that $(V_{j,k})_{jj} \to \infty$. The average amount of time spent in the $j$th
infected compartment is $\frac{1}{\alpha_{1,j}}$, and thus as $\alpha_{1,j} \to \infty$ the average
time spent in that compartment will go to zero.  This gives us an intuition for the idea that removing a compartment
from a system can be achieved through taking a limit.

It is easily observed that from Equations (6)-(8) of \cite{johnsonPLOSNTD} that the 
system with the $j$th infected compartment removed from the relapses of the first species is related to the whole system
through the Jacobians:
\[
	F_{j-1,j} = (F_{j,j})_{[j;j]} \mbox{ and } V_{j-1,j} = (V_{j,j})_{[j;j]}
\]
We can now apply the results of Corollary \ref{conthe} to find that
\begin{multline*}
	R_{0}^{j-1,j} = \rho(F_{j-1,j}(V_{j-1,j})^{-1})
	=\rho( (F_{j,j})_{[j;j]} (V_{j,j})_{[j;j]}^{-1}
	\\= \lim_{v_{jj} \to \infty} \rho (F_{j,j}(V_{j,j})^{-1})
	 = \lim_{\alpha_{1,j} \to \infty} R_{0}^{j,j}
\end{multline*}
Using \eqref{coupledrnot} we find that
\[
	R_{0}^{j-1,j} = \sqrt{\lim_{a_{1,j}\to\infty}(R_{0,1,j})^2 + (R_{0,2,j})^2}
\]
Apply \eqref{uncoupledrnot} and observe that
\[
	\lim_{a_{1,j}\to\infty}(R_{0,1,j})^2  = R_{0,1,j-1}
\]
Thus
\[
R_{0}^{j-1,j} = \sqrt{(R_{0,1,j-1})^2 + (R_{0,2,j})^2} 
\]
We can repeat this process, making the same observations and applying Corollary \ref{conthe}.  Thus we can say that 
when the first species undergoes $k-1<j-1$ relapses we get
\[
R_{0}^{k,j} = \sqrt{(R_{0,1,k})^2 + (R_{0,2,j})^2} 
\]
\section{Discussion}
	We have related the inverse of a matrix when the $i$th row and column are removed to the inverse of the 
	original matrix through a limit.  Such updating results generally have numerical uses, but in our case the existence
	of a limit is a complication.  Even the use of this method for giving and approximation to the updated matrix
	is impractical and inefficient, since it requires computation of the inverse of a larger matrix before the limit is
	taken.  However, we demonstrated an
	analytical application that allowed us to extend results to coupled systems of host and vectors in the spread of 
	a relapsing disease.  In particular, we were able to quantify how the two species undergoing a different number of 
	relapses affects the fundamental reproductive ratio for the disease.

\bibliographystyle{siamplain}
\bibliography{currentresearch}
\end{document}